\documentclass[onecolumn, 11pt]{IEEEtran}

\usepackage{amsmath}
\usepackage{amsthm}
\usepackage{amssymb}
\usepackage{subcaption}
\usepackage{mathtools}
\usepackage{cite}
\usepackage{color}
\usepackage{cite}
\usepackage{graphicx}

\newtheorem{theorem}{Theorem}

\newtheorem{example}{Example}
\newtheorem{assumption}{Assumption}

\newtheorem{definition}{Definition}
\newtheorem{lemma}{Lemma}
\newtheorem{corollary}{Corollary}

\renewcommand{\P}{\mathbb{P}}

\newcommand{\diag}{\textup{diag}}

   \def\vx{{\bf x}}

\def \e{\varepsilon}
\def \R{\mathbb{R}}
\def \dx{\,dx}
\def \E{\mathbb{E}}
\def \ones{\textbf{1}}

\def \calS{\mathcal{S}}
\def \calD{\mathcal{D}}

\def \R{\mathbb{R}}
\def \calN{\mathcal{N}}
\def \grad{\nabla}

\def \f{f} 
\def \F{F} 

\def \calF{\mathcal{F}}

\mathtoolsset{showonlyrefs}

\title{Distributed Gradient Methods for\\ Nonconvex Optimization:\\ Local and Global Convergence Guarantees}

\author{Brian Swenson,$^\dagger$ Soummya Kar,$^\ddagger$ H. Vincent Poor,$^\dagger$\\ Jos{\'e} M. F. Moura,$^\ddagger$ and Aaron Jaech$^{\star}$\thanks{\noindent
The work of B. Swenson and H. V. Poor was partially supported by the
Air Force Office of Scientific Research under MURI Grant FA9550-18-1-0502.\newline
$^\dagger$Department of Electrical Engineering, Princeton University, Princeton, NJ (\{bswenson, poor\}@princeton.edu),\newline
$^\ddagger$Department of Electrical and Computer Engineering, Carnegie Mellon University, Pittsburgh, PA (soummyak@andrew.cmu.edu and moura@ece.cmu.edu)\newline
$^\star$Facebook AI, Seattle, WA (ajaech@facebook.com)}
}

\begin{document}

\maketitle


\begin{abstract}
The article discusses distributed gradient-descent algorithms for computing local and global minima in nonconvex optimization. For local optimization, we focus on distributed stochastic gradient descent (D-SGD)---a simple network-based variant of classical SGD. We discuss local minima convergence guarantees and explore the simple but critical role of the stable-manifold theorem in analyzing saddle-point avoidance. For global optimization, we discuss annealing-based methods in which slowly decaying noise is added to D-SGD. Conditions are discussed under which convergence to global minima is guaranteed. Numerical examples illustrate the key concepts in the paper.
\end{abstract}

\vspace{-0em}
\section{Introduction}
Nonconvex optimization problems are prevalent in applications throughout control, signal processing, and machine learning. Modern applications involve unprecedented quantities of data generated by
a multitude of interconnected devices including mobile phones, IoT devices, self-driving vehicles, and networked cyber-physical systems. Due to limitations of the communication infrastructure it can be infeasible to collect the enormous amount of data generated by these devices to a centralized location for processing. 

This motivates the study of distributed (network-based) algorithms for nonconvex optimization. In a distributed optimization algorithm, a group of nodes (or agents) may communicate over an overlaid communication graph. Agents collaboratively optimize some collective function without any centralized coordination or data aggregation.\footnote{Note that this is distinct from some common decentralization schemes in which the computation is coordinated by a central node.} Distributed algorithms are flexible, robust, and efficient \cite{dimakis2010gossip,bullo2009distributed,nedic2009distributed,lian2017can}. They play a critical role in  IoT, edge computing, and sensor network applications, as well as large-scale parallel computing \cite{lian2017can}. 


While there are a wide variety of optimization techniques available, in this article we will focus on distributed implementations of simple first-order algorithms; namely, distributed gradient descent and variants thereof. 
Despite their simplicity, first-order methods are widely used in practice due to their ease of implementation, analytic tractability, and effectiveness in practical large-scale problems. 

We will consider the following setup: A group of  $N$ agents may communicate via an overlaid communication network. Each agent possesses some local objective function $f_n:\R^d\to\R$. It is desired to optimize the sum-function
\begin{equation} \label{eq:f-distributed}
F(x) := \sum_{n=1}^N f_n(x).
\end{equation}
Many problems of interest fall within this framework \cite{di2016next,rabbat2004distributed,lian2017can}. For example, in the context of machine learning, 
suppose that $\ell_n(x,y_{n})$ denotes the loss at agent $n$ given the parameter $x$ and datum $y_{n}$. Assuming the data at agent $n$ is distributed according to some distribution $\calD_{n}$, $f_n$ may correspond to the expected loss at agent $n$, i.e., $f_n(x) = \E_{y\sim\calD_n}(\ell_n(x,y))$, while \eqref{eq:f-distributed} corresponds to the expected global loss across all agents. 



Ideally, one would like to compute a global minimum of \eqref{eq:f-distributed}. However, when $F$ is nonconvex, computation of global minima can be challenging, and, in some important applications of interest, local minima are known to perform nearly as well as global minima \cite{ge2015escaping,kawaguchi2016deep,sun2015nonconvex,ge2016matrix}. Thus, it is prudent to begin by focusing on the problem of computing local minima of \eqref{eq:f-distributed}. 

Intuitively, gradient-based algorithms descend the objective function until they reach a point where the gradient is zero. That is to say, the set of limit points of gradient-based algorithms consists of the set of critical points. 
This set includes, of course, local minima, local maxima, and saddle points. Clearly, local maxima are not desirable limit points for an optimization algorithm, and it is easy to show that gradient-descent based methods do not actually converge to these points. However, gradient-based methods \emph{can} converge to saddle points, which can be quite problematic. 
Saddle points can be highly suboptimal and, in problems of interest, the number of saddle points can proliferate exponentially relative to the number of local minima as the problem size scales \cite{dauphin2014identifying}. 

The focus of this article will be on reviewing \emph{refined} convergence guarantees for distributed gradient algorithms. Specifically, we will consider distributed stochastic gradient descent (D-SGD) and focus on the following two fundamental questions:
\begin{itemize}
  \item [1.] Under what conditions can D-SGD be guaranteed to converge to local minima (or not converge to saddle points)?
  \item [2.] Can simple variants of D-SGD converge to global minima?
\end{itemize}

The fundamental issues involved in saddle point avoidance are most easily understood by considering continuous-time gradient flows.
In classical centralized gradient flow (GF), the key to understanding saddle-point behavior lies in the so-called \emph{stable-manifold theorem} \cite{coddington1955theory,shub2013global}.
Simply put, the main idea of the stable-manifold theorem is the following: Given any saddle point, there exists an associated \emph{stable manifold}---a special low-dimensional surface from which GF converges to the saddle point. In particular, GF converges to the saddle point \emph{if and only if} initialized on this surface. The (classical) stable-manifold theorem 
provides the essential structural information that enables analysis of saddle-point nonconvergence 
in centralized settings \cite{lee2016gradient,jin2017escape,murray2019revisiting}.

We will begin addressing the first question above by considering stable manifolds for \emph{distributed} gradient flow (DGF). While the stable manifold for DGF differs in some structural respects from the GF stable manifold, it serves the same essential role in facilitating analysis of saddle-point behavior in the distributed setting. In particular, we will see that DGF can only reach saddle points from a zero measure set of initializations.
 

After considering the stable-manifold theorem for DGF, we will turn our attention to discrete-time algorithms. We will consider distributed stochastic gradient descent (D-SGD)---a simple distributed variant of classical SGD. Under mild assumptions, D-SGD avoids saddle points and converges to local minima with probability 1. We will discuss conditions under which saddle points are avoided and highlight the key role of the DGF stable-manifold theorem in analyzing the saddle-point behavior of D-SGD.

Next, we turn to the problem of computing \emph{global} minima of \eqref{eq:f-distributed}.
Of course, the gradient of a function describes only a local property of the function and (stochastic) gradient descent is only capable of reliably locating local minima. In order to find global minima we must resort to other techniques. Simulated annealing (SA) is a popular method for locating global minima inspired by the annealing process in metallurgy wherein a metal is heated and then slowly cooled in order to freeze it into a low energy lattice configuration.
Analogously, the premise of SA algorithms is that slowly decaying noise may be added to a local search algorithm in order to escape local minima and seek out global minima. The slow reduction in added noise strength in SA corresponds to a slow ``cooling'' of the algorithm.
Simulated annealing algorithms were originally introduced for discrete combinatorial optimization \cite{kirkpatrick1983optimization}, with later variants being developed for continuous optimization in $\R^d$ \cite{gelfand1991recursive,yin1999rates,kushner1987asymptotic,chiang1987diffusion}.

In order to address the second question above, we will discuss gradient-descent based annealing algorithms for global optimization in $\mathbb{R}^d$. In particular, we will focus on a very simple generalization of D-SGD mentioned earlier in which slowly decaying (annealing) noise is added to the algorithm at each iteration. The algorithm, referred to as D-SGD + annealing, converges in probability to the set of global minima of $F$. 
We remark that the annealing algorithms we consider in this article are in the spirit of \cite{gelfand1991recursive} and are closely related to stochastic gradient Langevin dynamics \cite{zhang2017hitting,raginsky2017non,chiang1987diffusion,chen2019} currently popular in machine learning applications.


\vspace{.5em}
\noindent \textbf{Organization}. The remainder of the article is organized as follows. In Section \ref{sec:GD-central}, we begin by considering classical \emph{centralized} gradient algorithms. We review classical results for GF and SGD. 
In Section \ref{sec:DGF} we discuss DGF; we review convergence guarantees and basic structural properties of DGF. 
In Section \ref{sec:D-SGD} we consider D-SGD; we review guarantees for avoiding saddle points and converging to local minima. 
In Section \ref{sec:annealing-algorithms} we turn to the problem of locating global minima of \eqref{eq:f-distributed} using annealing methods. We review classical (centralized) results for gradient based annealing algorithms in Section \ref{sec:GD+annealing} and review D-SGD + annealing in Section \ref{sec:DGD+annealing}. In Section \ref{sec:sims} we consider numerical examples illustrating the key concepts in the paper. Section \ref{sec:conclusion} concludes the paper.

\vspace{.5em}
\noindent \textbf{Notation}. Throughout the paper, $\|\cdot\|$ denotes the standard Euclidean norm, $\textup{dist}(x,y) = \|x-y\|$ is the distance between points $x$ and $y$, and $\textup{dist}(x,S) = \inf_{y\in S}\textup{dist}(x,y)$ gives the distance between a point $x$ and set $S$. Given a random process $\{y_k\}_{k\geq 1}$ taking the general recursive form $y_{k+1} = G(y_k, \xi_k)$, where $\xi_k$ is a random variable, $k$ is the iteration, and $G(\cdot,\cdot):\R^d\times \R^d\to\R^d$ is the iteration map, we let $\calF_k = \sigma(\xi_{1},\ldots,\xi_{k-1})$ denote the $\sigma$-algebra representing the information available at time $k$ and let $\E(\cdot\vert\calF_k)$ represent the associated conditional expectation.\footnote{Informally, $\calF_k$ contains information about the outcome of all random events that have occurred up to time $k$ and the notation $\E(\cdot\vert\calF_k)$ represents a conditional expectation taken with respect to the random variables inside $\sigma(\cdot)$ in the definition of $\calF_k$.}

\section{Centralized Gradient Algorithms: Foundations and Intuition} \label{sec:GD-central}
Before studying properties of distributed algorithms in nonconvex optimization, it will be helpful to first review properties of centralized gradient algorithms. 
We begin by considering classical continuous-time GF. 

\subsection{Centralized Gradient Flow} \label{sec:GF}
Given a differentiable function $\f:\R^d\to\R$, the well-known stochastic gradient descent algorithm is given by the recursion 
\begin{equation} \label{eq:SGD-centralized}
x_{k+1} = x_k - \alpha_k(\nabla f(x_k) + \xi_{k}),
\end{equation}
where $\xi_k$ is unbiased noise, i.e., a random variable satisfying $\mathbb{E}(\xi_k\vert\calF_k) = 0$, 
$\alpha_k$ denotes the step size (or learning rate), and $\calF_k$ represents the information available at iteration $k$. Taking an expectation on both sides and rearranging terms we have
\begin{equation} \label{eq:GD-central2}
\frac{\E(x_{k+1} - x_k\vert \calF_k)}{\alpha_k} = -\nabla \f(x_k).
\end{equation}

\noindent If the step size $\alpha_k$ is taken to zero at an appropriate rate (see Assumption \ref{a:step-size1} below), then \eqref{eq:GD-central2} represents an Euler discretization of the \emph{gradient flow} (GF) differential equation, given by
\begin{equation} \label{eq:GF-centralized}
\dot \vx(t) = -\nabla \f(\vx(t)).
\end{equation}
A differentiable function $\vx:[0,\infty)\to\R^d$ is said to be a solution to \eqref{eq:GF-centralized} with initial condition $x_0\in\R^d$ if $\vx(0) = x_0$ and $\vx(t)$ satisfies \eqref{eq:GF-centralized} for all $t>0$. We will use the convention of representing solutions of ODEs with bold face text. 
The relationship between \eqref{eq:SGD-centralized} and \eqref{eq:GF-centralized} is made rigorous using tools from the field of  stochastic approximation theory \cite{benaim1999dynamics}, 
which deals with analysis techniques for studying stochastic discrete-time systems by considering the corresponding (deterministic) continuous time ODE.
The advantage of studying the ODE is that fundamental properties of the system can be much easier to understand. This has motivated a recent trend in ODE-based methods for algorithm analysis and development \cite{krichene2015accelerated,su2014differential,wibisono2016variational}. 

In this section we will review the basic properties of \eqref{eq:GF-centralized} in nonconvex optimization.
The following assumption gives a basic condition under which GF is well defined \cite{coddington1955theory}. We use the notation $f\in C^k$ for integer $k\geq 0$ to indicate that $f$ is $k$-times continuously differentiable. 
\begin{assumption} \label{a:lip-grad}
$\f$ is $C^1$ and has Lipschitz continuous gradient, i.e., there exists a constant $K>0$ such that $\|\nabla \f(x) - \nabla \f(y)\| \leq K\|x-y\|$ for all $x,y\in\R^d$.
\end{assumption}
In order for GF to be useful, we would like it to converge to some point or set. The following assumption ensures that GF does not flow outward indefinitely. The assumption simply asserts that, asymptotically, the negative gradient points inwards. 
\begin{assumption} \label{a:coercive}
There exists a radius $R>0$ and constant $C>0$ such that $
\big\langle  \frac{\nabla f(x)}{\|\nabla f(x)\|},\frac{x}{\|x\|} \big\rangle \geq C $
for all $\|x\|\geq R$.
\end{assumption}

Note that the only points at which GF may rest are points where the right hand size of \eqref{eq:GF-centralized} is zero. These are precisely the set of critical points of $f$. The following standard result shows that GF converges to this set.
\begin{theorem} [Convergence to Critical Points] \label{thrm:GF-CP} Suppose that $f$ satisfies Assumption \ref{a:lip-grad} and $\vx$ satisfies \eqref{eq:GF-centralized}. Then every limit point of $\vx$ is a critical point of $f$. If Assumption \ref{a:coercive} also holds, then $\vx$ converges to the set of critical points of $f$.
\end{theorem}

Of course, the set of critical points of $f$ consists of local maxima, local minima, and saddle points. Thus, some limit points permitted under Theorem \ref{thrm:GF-CP} can be extremely suboptimal. We are interested in refinements of Theorem \ref{thrm:GF-CP} that show that the GF typically converges to local minima.
It is not difficult to show that GF does not converge to local maxima.  However, saddle points require a subtler treatment.
In order to build intuition, consider the following simple example.
\begin{example} \label{example1}
Let $\f:\R^2\to\R$ be given by the quadratic function $\f(x) = \frac{1}{2}(x_1^2 - x_2^2).$
Note that $\f$ has a saddle point, which we will denote by $x^*$, at the origin. A plot of this function is shown in Figure \ref{fig:example1-saddle}, and a plot of the gradient-descent vector field is shown in Figure \ref{fig:example1-grad-vec-field}.

\begin{figure}[h]
  \begin{subfigure}[h]{0.35\textwidth}
    \includegraphics[width=1\textwidth]{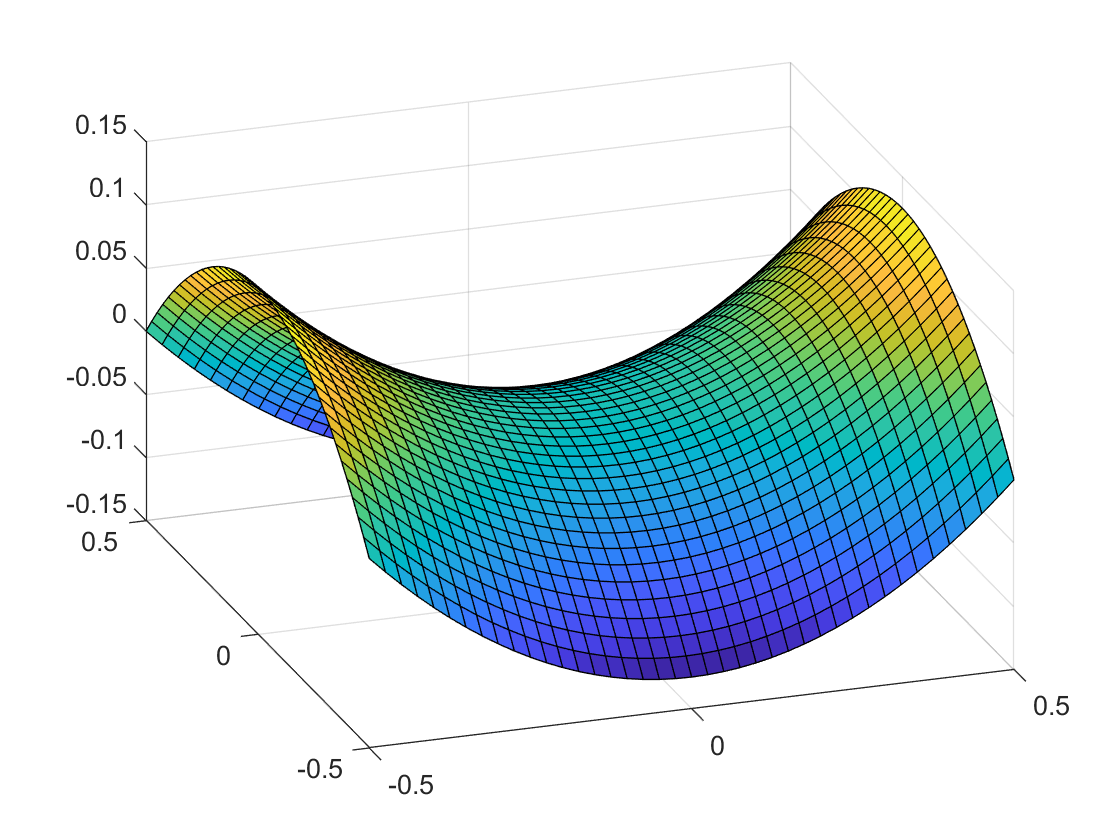}
    \caption{}
    \label{fig:example1-saddle}
  \end{subfigure}
  \begin{subfigure}[h]{0.31\textwidth}
    \includegraphics[width=1\textwidth]{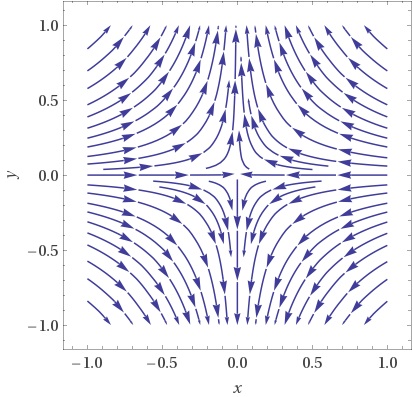}
    \caption{}
    \label{fig:example1-grad-vec-field}
  \end{subfigure}
  \begin{subfigure}[h]{.31\textwidth}
    \includegraphics[width=1\textwidth]{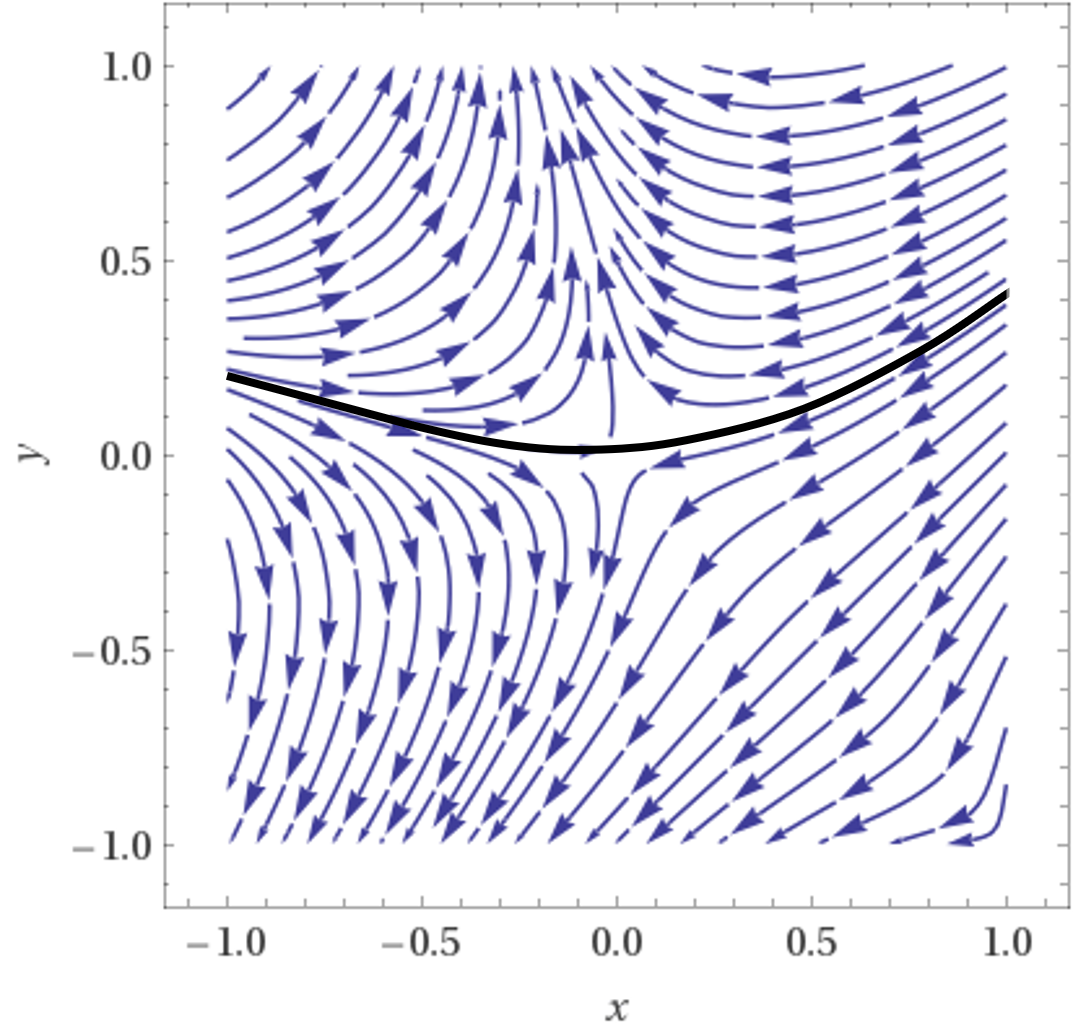}
    \caption{}
    \label{fig:example3}
  \end{subfigure}
\end{figure}

We would like to characterize the set of initial conditions from which solutions to \eqref{eq:GF-centralized} converge to the saddle point $x^*$. Formally, this set is given by $\{x_0\in\R^2: \vx(0) = x_0 \mbox{ and } \vx \mbox{ satisfies } \eqref{eq:GF-centralized} \implies \vx(t)\to x^*\}$,
and (for reasons soon to become clear) we will refer to it as the \emph{stable manifold} associated with $x^*$.

Since $f$ is quadratic, the GF system associated with $f$ is given by the linear system $
\dot \vx(t) = -A\vx(t)$ where $A=\diag(1,-1)$ is the diagonal matrix with $1$ and $-1$ on the diagonal.
Using basic tools from linear systems theory, we see the solution to this system given an initial condition $\vx(0) = x_0$ is given by
$$
\vx(t) = e^{-At}x_0 =
\begin{pmatrix}
e^{-t} & 0 \\
0 & e^{t}
\end{pmatrix}x_0.
$$
Note that if we initialize a GF trajectory on the $x_1$ axis, i.e., with $\vx(t_0) = x_0$, $x_{0}\in \{x\in\R^2:x_2=0\}$, then $\vx(t) \to x^*$.  However, if we initialize $x_0\not\in\{x:x_2=0\}$, then not only does $\vx(t)\not\to x^*$, but we have $|\vx_2(t)|\to \infty$ as $k\to\infty$.
\end{example}

The next example is an extension of Example \ref{example1} that illustrates how the stable manifold depends on the dimension of the problem and the structure of the saddle point.
Before presenting the example, we recall the following standard notation. For an integer $p\geq 1$, $\ones_{p}$ is the $p$-dimensional vector of all ones. Given vectors $a\in \R^{d_1}$ and $b\in \R^{d_2}$, $\diag(a,b)$ gives the $(d_1+d_2)\times(d_1+d_2)$ diagonal matrix with the elements of $a$ and $b$ on the diagonal.
\begin{example}
Let $d\geq 2$ and $q\in \{1,\ldots,d\}$.
Suppose $f:\R^d\to \R$ is given by $f(x) = \frac{1}{2}x^T A x$ where $A = \diag(\ones_{d-q},-\ones_{q})$
so that $q$ denotes the number of \emph{negative} eigenvalues of $\nabla^2 f(0) = A$. Because $f$ is quadratic, \eqref{eq:GF-centralized} is given by the linear system $\dot \vx = -A\vx.$  Since $A$ is diagonal, the solution to this system with initial condition $x_0\in \R^d$ is simply given by
$$
\vx(t) = \begin{pmatrix}
\ones_{d-q}e^{-t}\\
\ones_{q}e^{t}
\end{pmatrix}x_0.
$$
Given this form of the solution, it is clear that the stable manifold for $x^*$ is the set $\calS = \{x\in \R^d: x_i=0,~i=d-q+1,\ldots,d\}.$
Note that $\calS$ is the subspace spanned by all eigenvectors of $A$ with positive (or, more generally, nonnegative) eigenvalue.

We are interested in understanding when the stable manifold will be a ``small'' set in some sense. One way to gauge the size of $\calS$ is by its dimension. Recalling that $q$ is the number of negative eigenvalues and that $\calS$ is the span of all eigenvectors with nonnegative eigenvalues, we have $\dim \calS = d-q$.
\end{example}

The intuition garnered from these elementary examples will generalize to a wide range of saddle points. We highlight the following key observations from these examples: (i) While it is possible to converge to the saddle point $x^*$, this can only be accomplished by initializing on a special low-dimensional surface (referred to as the stable manifold of $x^*$). (ii) The stable manifold is an invariant set. That is, if a trajectory $\vx(t)$ lies in the stable manifold at some time $t_0\geq 0$, then $\vx(t)$ remains in the stable manifold for all $t\geq t_0$. (iii) The stable manifold is a ``repulsive'' set. If a trajectory $\vx(t)$ is not initialized \emph{precisely} on the stable manifold (e.g., if $\vx_2(0) = \e > 0$ in Example \ref{example1}), then it is pushed away from the stable manifold.
(iv) In the above examples, the GF system was a linear system and the stable manifold was a linear subspace of $\R^d$. This occurred because we chose $f$ to be quadratic. In general, when $f$ is not quadratic, the stable manifold will be a smooth nonlinear surface.
(v) The stable manifold has dimension $d-q$, where $q$ is the number of negative eigenvalues in $\nabla^2 f(x^*)$. This relationship was observed when $\calS$ was a linear surface, but will hold generally when $\calS$ is an arbitrary smooth surface.

These observations generalize, and typify the behavior of GF near a broad class of saddle points known as \emph{regular} saddle points.\footnote{In the optimization literature, saddle points satisfying this condition are often referred to as \emph{nondegenerate}. Because we will discuss nonconvergence to these points, we prefer the term regular to avoid double negatives.} 
\begin{definition}[Regular Saddle Point] \label{def:strict-saddle}
 We say that a saddle point $x^*$ of $f$ is \emph{regular} (or \emph{nondegenerate}) if $\nabla^2 f(x^*)$ exists and is nonsingular.
\end{definition}
Regular saddle points serve two purposes: (1) They ensure that $\nabla^2 f(x^*)$ has at least one negative eigenvalue, which implies that $f$ has at least one genuine descent direction at $x^*$. (2) They greatly simplify the analysis.
\footnote{Some authors have considered the notion of a \emph{strict} saddle point, i.e., a critical point $x^*$ of $f$ where $\nabla^2 f(x^*)$ exists and has at least one negative eigenvalue \cite{jin2017escape}. The notion of a strict saddle is weaker than a regular saddle (although in instances where the function $f$ is generated by random data the distinction between these two notions is typically one of probability zero). Stable manifold results exist for strict saddles in centralized settings, but not in distributed, hence we will focus on regular saddle points here.}

The following theorem generalizes the intuition from the above examples and shows that a stable manifold exists near regular saddle points. This result refines Theorem \ref{thrm:GF-CP} by showing that the critical point reached by GF cannot generically be a regular saddle point. A proof of the theorem can be found in \cite{coddington1955theory}.
\begin{theorem}[Stable-Manifold Theorem for GF] \label{thrm:stable-manifold-CT-central}
Suppose that $f$ is $C^2$ and satisfies Assumption \ref{a:lip-grad} and that $x^*$ is a regular saddle point of $f$. Let $q$ denote the number of negative eigenvalues of $\nabla^2 f(x^*)$. Then there exists a smooth $(d-q)$-dimensional surface $\calS$  such that a solution $\vx$ to \eqref{eq:GF-centralized} can only converge to $x^*$ if it is initialized on $\calS$.
\end{theorem}

It is important to note that, in the above theorem, since $x^*$ is a regular saddle, we have $q\geq 1$ and so $\calS$ has dimension at most $d-1$. Thus, $\calS$ is a Lebesgue measure zero set in $\R^d$.
An alternate, but intuitively satisfying way to restate the above result is that, if GF is randomly initialized according to a ``nice'' distribution,\footnote{Meaning a distribution whose underlying probability measure is absolutely continuous with respect to the Lebesgue measure.}
then GF avoids regular saddle points with probability 1.

To illustrate the structure of the stable manifold near a function that is not quadratic but possesses a regular saddle point, suppose $f(x) = x_1^2 - x_2^2 + x_1^2x_2 + x_1x_2^2$, and note that $f$ has a regular saddle point at the origin. The GF vector field for this function is displayed in Figure \ref{fig:example3}. The stable manifold is illustrated in black.

\subsection{Stochastic Gradient Descent} \label{sec:SGD-central}
Stochastic gradient descent is given by the recursive process \eqref{eq:SGD-centralized}.
The following is a standard assumption for the noise in SGD.
\begin{assumption}\label{a:noise-zero-mean}
$\E(\xi_k\vert\calF_k) = 0$ and $\E\left(\|\xi_k\|^2\vert \calF_k\right) \leq C$
for some $C>0$ and all $k\geq 1$.
\end{assumption}
This setup is quite broad and includes, among other things, mini-batch variants of SGD common in machine learning and neural network training as well as SGD for online learning \cite{shalev2014understanding}. 
Due to noise, SGD can only be guaranteed to converge if a decaying step size is used. We will assume that the step size takes the following form.
\begin{assumption}\label{a:step-size1}
$\alpha_k = \Theta(k^{\tau_\alpha})$, $\tau_\alpha \in (1/2,1]$.
\end{assumption}
In the above assumption, we use the asymptotic notation $\alpha_k = \Theta(k^{\tau_\alpha})$ to mean that for some constants $c_1,c_2>0$, $ c_1k^{\tau_\alpha}\leq \alpha_k \leq c_2 k^{\tau_\alpha}$ for all $k$ sufficiently large.
Step sizes in this range ensure that $\alpha_k$ decays quickly enough to average out noise, but slowly enough so that $x_k$ can adequately explore the state space.

The following technical assumption is standard for stochastic approximation algorithms \cite{davis2018stochastic}. We recall that a set is said $S\subset \R$ to be dense in $\R$ if for every point $x\in\R$ one may select a point $y\in S$ arbitrarily close to $x$. 
\begin{assumption} \label{a:CP-meas-zero}
Let $\text{CP}_f\subset \R^d$ denote the set of critical points of $f$. The set $\R\backslash f(\text{CP}_f)$ is dense in $\R$.
\end{assumption}

This assumption is quite mild and typically satisfied by functions encountered in practice \cite{davis2018stochastic}. 
The following theorem establishes the convergence of SGD to critical points. 
\begin{theorem} \label{thrm:SGD-conv-to-CP}
Suppose that Assumptions \ref{a:lip-grad}--\ref{a:CP-meas-zero} are satisfied. Then SGD converges to the set of critical points of $f$ with probability 1.
\end{theorem}

Note that Theorem~\ref{thrm:SGD-conv-to-CP} extends the convergence result of Theorem~\ref{thrm:GF-CP} from GF to SGD making the additional Assumptions \ref{a:noise-zero-mean}--\ref{a:CP-meas-zero}.
In the same spirit, the result follows using the ODE method of stochastic approximation \cite{benaim1999dynamics,davis2018stochastic}, which allows one to use GF as a surrogate for studying SGD.\footnote{The result follows directly from Theorem 4.2 in \cite{davis2018stochastic} where it may be verified that iterate sequence is bounded under Assumptions \ref{a:lip-grad}--\ref{a:step-size1}.}
Of course, convergence to critical points is rudimentary convergence criterion. We are interested in refinements of this result, and, in particular, understanding the (non)convergence of SGD to saddle points. 

In order to ensure nonconvergence to saddle points, we must make an additional assumption about the noise process. Informally, the following assumption states that the random variable $\xi_k$ in \eqref{eq:SGD-centralized} perturbs in all directions. The assumption ensures that the noise in SGD will knock the process off of any ``bad'' low-dimensional sets that could lead to a saddle point (i.e., a stable manifold). In the following assumption, we use the notation $(a)^+ := \max\{a,0\}$.
\begin{assumption}[Minimum Excitation] \label{a:pemantle-noise}
  The noise process satisfies $\E( (\xi_k^T\theta)^+\vert \calF_k) \geq C$ for some $C>0$ and every unit vector $\theta$.
\end{assumption}

We emphasize that this assumption is easily satisfied, for example, by any noise with positive definite covariance. 
The following theorem establishes that SGD does not converge to regular saddle points.
\begin{theorem} \label{thrm:SGD-nonconvergence}
Suppose that $f$ is $C^3$, that $x^*$ is a regular saddle point of $f$, and that Assumptions \ref{a:lip-grad}--\ref{a:pemantle-noise} are satisfied. Then $\P(x_k \to x^*) = 0.$
\end{theorem}

Note that Theorem \ref{thrm:SGD-nonconvergence} refines Theorem \ref{thrm:SGD-conv-to-CP}. The refinement is obtained by adding the minimum excitation assumption and stipulating that $f$ be smoother than previously assumed.
Theorem \ref{thrm:SGD-nonconvergence} is proved in \cite{pemantle1990nonconv}.\footnote{Reference \cite{pemantle1990nonconv} proves this result assuming bounded noise. However, it is straightforward to extend the arguments to handle noise with bounded variance, e.g., \cite{swenson2020saddles}.} The proof again relies on studying the underlying GF ODE. 
The complete proof of this result is fairly involved; however, the basic intuition underlying the proof is simple to grasp and elucidates the role of saddle-point structure and the GF stable-manifold theorem. We will briefly discuss this intuition now. We remark that the intuition underlying saddle-point nonconvergence of the \emph{distributed} variant of SGD (D-SGD), to be discussed in Section \ref{sec:D-SGD}, will be similar.

The main idea is that Theorem \ref{thrm:SGD-nonconvergence} follows as a consequence of the stable-manifold theorem for GF (Theorem \ref{thrm:stable-manifold-CT-central}). This follows from two key observations.
First, note that the \emph{mean} update step for SGD satisfies \eqref{eq:GD-central2}.
Thus, the mean update step of SGD is a discretization of \eqref{eq:GF-centralized} with decaying step size. Consequently, \eqref{eq:GF-centralized} is the asymptotic mean field for \eqref{eq:SGD-centralized}, and the asymptotic behavior of \eqref{eq:SGD-centralized} is determined by properties of \eqref{eq:GF-centralized}.
The second key observation is that the stable manifold for GF is a Lyapunov \emph{unstable} set. That is, if a trajectory is not initialized precisely on the stable manifold, then it is pushed away from the stable manifold. This property was illustrated in the Example \ref{example1} and holds generally for the stable manifold near any regular saddle.

Intuitively then, the GF stable manifold is precisely the object from which solutions to \eqref{eq:SGD-centralized} are repelled near $x^*$. Letting $\calS$ denote the stable manifold established in Theorem \ref{thrm:stable-manifold-CT-central}, Theorem \ref{thrm:SGD-nonconvergence} follows by showing that (i) noise satisfying Assumption \ref{a:pemantle-noise} pushes $x_k$ off of and away from $\calS$, then (ii) the instability of $\calS$ eventually forces $x_k$ away from $\calS$ forever. See \cite{pemantle1990nonconv} for more details.

\section{Distributed Gradient Flow} \label{sec:DGF}

We now consider distributed gradient processes for locating local minima. Analogous to the centralized setting, it is advantageous to begin by considering \emph{continuous-time} DGF. 

\subsection{DGF} \label{sec:DGF-intro}

Let $\vx_n(t)$ denote agent $n$'s estimate of the solution to \eqref{eq:f-distributed} at time $t$.
DGF is given by the differential equation
\begin{equation} \label{eq:DGF}
\dot \vx_n(t) = \beta_t \sum_{\ell \in \Omega_n} (\vx_\ell(t) -\vx_n(t)) - \alpha_t \grad f_n(\vx_n(t)).
\end{equation}
%
The right hand side of \eqref{eq:DGF} consists of two terms: A consensus term $\beta_t \sum_{\ell \in \Omega_n} (\vx_\ell(t) -\vx_n(t))$, and a local gradient-descent term $-\alpha_t \grad f_n(\vx_n(t))$.
As suggested by the name, the consensus term encourages agents' state estimates to tend towards a common value as $t\to\infty$.
(For example, if we suppose that $f_n \equiv 0$ for all $n$ then \eqref{eq:DGF} reduces to classical consensus dynamics \cite{olfati2007consensus,dimakis2010gossip}.) The gradient term encourages each agent to descent the gradient of their local objective $f_n$.
In order to ensure that solutions to \eqref{eq:DGF} are well-defined, we will assume that each $f_n$ satisfies Assumption \ref{a:lip-grad} with $f= f_n$.

\subsection{Consensus and Convergence to Critical Points} \label{sec:DGF-CPs}
In order to guarantee the convergence of \eqref{eq:DGF} we must make a few additional assumptions pertaining to the distributed setup. First, we assume that the communication graph is connected so that information may disseminate freely between nodes.\footnote{For simplicity, we will restrict attention to undirected time-invariant graphs throughout the paper.}
\begin{assumption} \label{a:G-connected-undirected}
  The graph $G=(V,E)$ is undirected and connected.
\end{assumption}


We will assume that the weight sequences take the following form.
\begin{assumption} \label{a:step-size-CT}
  $\alpha_t = \Theta( t^{-\tau_\alpha})$ and $\beta_t = \Theta( t^{-\tau_\beta})$, $0\leq \tau_\beta < \tau_\alpha \leq 1$.
\end{assumption}

Under the above assumptions, agents achieve consensus and agents' states $\vx_n(t)$ converge to a critical point of $F$, as stated in the following theorem. A proof of the theorem can be found in \cite{smkp2020TAC}\cite{swenson2019allerton}.
\begin{theorem}
\label{thrm:cont-conv-cp}
Suppose $\{\vx_n(t)\}_{n=1}^N$ is a solution to \eqref{eq:DGF} with arbitrary initial condition. Suppose that Assumptions \ref{a:lip-grad}--\ref{a:coercive} hold with $f_n=f$ for each $n$, and Assumptions \ref{a:G-connected-undirected}--\ref{a:step-size-CT} hold. Let $F$ be given by \eqref{eq:f-distributed}. Then for each agent $n$ we have: (i) $\lim_{t\to\infty} \|\vx_n(t) - \vx_\ell(t)\| = 0$, for all $\ell = 1,\ldots,N$, and (ii) $\vx_n(t)$ converges to the set of critical points of $F$.
\end{theorem}

Theorem \ref{thrm:cont-conv-cp} is an extension of Theorem \ref{thrm:GF-CP} to the distributed setting. Note that to obtain this extension, the theorem only requires the additional assumptions that the communication graph is connected and the consensus vs optimization weights $\beta_t$ and $\alpha_t$ are properly balanced. 

\subsection{Nonconvergence to Saddle Points} \label{sec:DGF-SPs}
We will now consider a refinement of Theorem \ref{thrm:cont-conv-cp} showing that, while it is possible to converge to saddle points, DGF typically avoids them. As in the case of GF, this can be accomplished by studying stable manifolds.

The stable manifold for DGF differs from the centralized case in a few important respects. First, in order to establish the existence of the stable manifold it is necessary to make the following mild technical assumption. In the statement of the assumption, $x^*$ refers to a saddle point of interest. 
\begin{assumption}[Continuity of Eigenvectors] \label{a:eigvec-continuity}
For each $n$, the eigenvectors of $\nabla^2 f_n(x)$ are continuous at $x^*$ in the sense that, for each $x$ near $x^*$, there exists an orthonormal matrix $U_n(x)$ that diagonalizes $f_n(x)$ such that $x\mapsto U_n(x)$ is continuous at $x^*$.
\end{assumption}
This assumption is relatively innocuous and should be satisfied by most functions encountered in practice; however, it is required to rule out certain pathological cases which can arise in the distributed setting (but not in the centralized) \cite{smkp2020TAC}.

The stable manifold for DGF also differs from the stable manifold for GF in that it has a time dependence.
To see this, note that the right-hand side of \eqref{eq:DGF} depends not only on the state of the process $\vx(t)$, but also on time via the time-varying weights $\alpha_t$ and $\beta_t$. As a result, the asymptotic behavior of a DGF trajectory with initialization $\vx_n(t_0) = x_{n,0}$, $n=1,\ldots,N$, 
depends not only on the initial state $x_0 = \{x_{n,0}\}_{n=1}^N$, but also on the initial time $t_0$. 
Consequently, the stable manifold for DGF will be dependent time-dependent 

The following theorem establishes the existence of a stable manifold near regular saddle points and characterizes the structure of the stable manifold. A proof of this result can be found in \cite{smkp2020TAC,swenson2019allerton}.
\begin{theorem} \label{thrm:non-convergence-DGF}
Suppose that the hypotheses of Theorem \ref{thrm:cont-conv-cp} hold. 
Suppose also that each $f_n$ is $C^2$, that $x^*$ is a regular saddle point of $\F$ satisfying Assumption \ref{a:eigvec-continuity}.
Let $q$ denote the number of negative eigenvalues of the Hessian $\nabla^2 F(x^*)$. Then for all $t_0$ sufficiently large there exist a smooth surface $\calS_{t_0}\subset \R^{Nd}$ with dimension $(Nd-q)$ such that the following holds: A solution $(\vx_n(t))_{n=1}^N$ to \eqref{eq:DGF} may converge to $x^*$ in the sense that $\vx_n(t)\to x^*$ for some $n$ (then every $n$), only if $(\vx_n(t))_{n=1}^N$ is initialized on $\calS_{t_0}$, i.e., $(\vx_n(t_0))_{n=1}^N = x_0\in\R^{Nd}$ with $x_0\in S_{t_0}$.
\end{theorem}

Theorem \ref{thrm:non-convergence-DGF} refines Theorem \ref{thrm:cont-conv-cp}, showing that the critical point obtained by DGF cannot generally be a regular saddle point. The refinement is obtained by making the additional assumptions that $f$ is $C^2$ and Assumption \ref{a:eigvec-continuity} holds near the saddle point.\footnote{The assumption that $f$ is $C^2$ need only hold locally near $x^*$ \cite{smkp2020TAC}.}
The following corollary illustrates a simple condition under which DGF converges to local minima. The corollary is an immediate consequence of Theorem \ref{thrm:non-convergence-DGF}. In the statement of the theorem, we remark that ``almost every'' is in the sense of the Lebesgue measure in $\R^{Nd}$.

\begin{corollary}
Suppose that the hypotheses of Theorem \ref{thrm:non-convergence-DGF} hold and that every saddle point of $\F$ is regular and satisfies Assumption \ref{a:eigvec-continuity}. Then, for any initial time $t_0$, $\vx(t)$ converges to the set of local minima from almost every initialization $x_0\in \R^{Nd}$.
\end{corollary}


\section{Distributed Stochastic Gradient Descent} \label{sec:D-SGD}
We will now consider D-SGD---a simple distributed variant of the classical SGD algorithm \eqref{eq:SGD-centralized}.
Before introducing the D-SGD algorithm, we will briefly review closely related work on this nascent topic. Reference \cite{daneshmand2018second} considers a variant of discrete-time DGD with constant step size and shows convergence to a neighborhood of second order stationary points for sufficiently small step sizes for almost all initial conditions. A critical point $x^*$ of $f$ is called second order stationary if $\nabla^2 f(x^*)$ is positive semidefinite.  Thus, for example, if all saddle points are regular per Definition \ref{def:strict-saddle}, then this implies convergence to local minima. Reference \cite{daneshmand2018second} also considers distributed gradient tracking algorithms and saddle-point nonconvergence. References \cite{vlaski2019distributed1,vlaski2019distributed2} consider a variant of D-SGD with constant step size and show polynomial-time escape from saddle points. Reference \cite{hong2018gradient} considers primal-dual based methods for distributed optimization and shows convergence to second order stationary points from almost all initializations.

\subsection{D-SGD Algorithm} \label{sec:D-SGD-intro}
For integers $k\geq 1$, let $x_n(k)\in \R^d$ denote the estimate that agent $n$ maintains of an optimizer of \eqref{eq:f-distributed} at iteration $k$. The D-SGD algorithm is given by the recursion
\begin{equation} \label{eq:S-DGD}
x_n(k+1) = x_n(k) + \beta_k\sum_{\ell\in \Omega_n} (x_\ell(k) - x_n(k)) - \alpha_k( \nabla f_n(x_n(k)) + \xi_n(k)),
\end{equation}
for each agent $n=1,\ldots,N$, where $\xi_n(k)$ denotes a random variable. The random variable $\xi_n(k)$ typically represents gradient measurement noise
but may also represent noise that is deliberately introduced to aid in escaping saddle points, e.g., \cite{jin2017escape}.
Similar to \eqref{eq:DGF}, the right hand side of \eqref{eq:S-DGD} consists of a consensus term $\beta_k\sum_{\ell\in \Omega_n} (x_\ell(k) - x_n(k))$ and a (stochastic) gradient-descent term $-\alpha_k( \nabla f_n(x_n(k)) + \xi_n(k))$.


\subsection{Consensus and Convergence to Critical Points} \label{sec:SGD-consensus}
In studying the convergence of D-SGD, we will retain Assumption \ref{a:G-connected-undirected} used for DGF
and (analogous to Assumption \ref{a:step-size-CT}) we will assume the discrete-time weight parameters take the following form. 
\begin{assumption} \label{a:step-size-DT}
$\alpha_k = \Theta\left( k^{-\tau_\alpha}\right)$ and $\beta_k = \Theta\left(k^{-\tau_\beta}\right)$ with $0\leq \tau_\beta < \tau_\alpha$, $\tau_\alpha\in (1/2, 1]$.
\end{assumption}



This assumption is similar to Assumption \ref{a:step-size1} from SGD, but here we also account for the consensus weight $\beta_k$.
The next theorem shows that D-SGD converges to critical points of $F$.
\begin{theorem}
\label{thrm:discrete-conv}
Let $\{\{x_n(k)\}_{n=1}^N\}_{k\geq 1}$ be a D-SGD process \eqref{eq:S-DGD}. 
Suppose that Assumptions \ref{a:lip-grad}--\ref{a:coercive} hold with $f_n = f$. Assume that for each $n$, $\xi_n(k)$ is independent from $\xi_\ell(k)$, $\ell\not= n$ and satisfies Assumption \ref{a:noise-zero-mean} with $\xi_n(k) = \xi_k$. Suppose that Assumption \ref{a:CP-meas-zero} holds with $F = f$. 
Suppose also that Assumptions \ref{a:G-connected-undirected} and \ref{a:step-size-DT} hold. Then with probability 1, for each agent $n$ the following holds: (i) Asymptotic consensus is achieved in the sense that $\lim_{k\to\infty} \|x_n(k) - x_\ell(k)\|=0$ for all $\ell=1,\ldots,N$, and (ii) $\{x_n(k)\}_{n=1}^N$ converges to the set of critical points of $F$.
\end{theorem}
Theorem \ref{thrm:discrete-conv} extends the convergence result of Theorem \ref{thrm:SGD-conv-to-CP} to the distributed setting. The extension requires the additional assumptions that the communication graph is connected and the consensus vs optimization weight parameters are balanced. The theorem is proved in \cite{swenson2020saddles} by using Theorem \ref{thrm:non-convergence-DGF} and ODE-based stochastic approximation techniques.

\subsection{Nonconvergence to Saddle Points} \label{sec:D-SGD-SPs}
In the case of (centralized) SGD, we saw that under a minimum excitation condition (Assumption \ref{a:pemantle-noise}), the gradient noise knocked the process away from the underlying stable manifold, allowing the algorithm to escape from saddle points with probability one.

The following result shows that D-SGD does not converge to saddle points under a similar minimum excitation assumption.
\begin{theorem}
\label{thrm:nonconvergence-DT-DSGD}
Suppose 
that the hypotheses of Theorem  \ref{thrm:discrete-conv} hold. Suppose also that $F\in C^3$, $x^*$ is a regular saddle point of $F$ satisfying Assumption \ref{a:eigvec-continuity}, and for each $n$, $\{\xi_n(k)\}_{k\geq 1}$ satisfies Assumption \ref{a:pemantle-noise} with $\xi_n(k) = \xi_k$.
Then, regardless of initialization, for each $n$, $\P(x_n(k)\to x^*) = 0.$
\end{theorem}
Note that Theorem \ref{thrm:nonconvergence-DT-DSGD} extends the result of Theorem \ref{thrm:SGD-nonconvergence} to the distributed setting and refines Theorem \ref{thrm:discrete-conv}. The extension of Theorem \ref{thrm:SGD-nonconvergence} is again obtained by making the additional assumptions that the communication graph is connected and the consensus and gradient weights $\beta_k$ and $\alpha_k$ are balanced. The refinement of Theorem \ref{thrm:discrete-conv} is obtained by making the additional assumptions that $F$ is $C^3$ and Assumptions \ref{a:pemantle-noise} and \ref{a:eigvec-continuity} hold. Note the analogous relationship between Theorems \ref{thrm:SGD-conv-to-CP} and \ref{thrm:SGD-nonconvergence} in the centralized setting. 

In the centralized setting, nonconvergence of SGD (Theorem \ref{thrm:SGD-nonconvergence}) was obtained from the continuous-time flow via Theorem \ref{thrm:stable-manifold-CT-central}. In the distributed setting, the relationship is analogous. Nonconvergence of D-SGD (Theorem \ref{thrm:nonconvergence-DT-DSGD}) is obtained from the continuous-time flow via Theorem \ref{thrm:non-convergence-DGF}.
The idea is fundamentally the same as the proof of Theorem \ref{thrm:SGD-nonconvergence} discussed in Section \ref{sec:SGD-central}.
Noise satisfying Assumption \ref{a:pemantle-noise} at each agent perturbs the D-SGD process off of and away from the stable manifold. The stable manifold is inherently unstable (because $x^*$ is a regular saddle point). This inherent instability eventually pushes the D-SGD process away from the saddle point forever. See \cite{swenson2020saddles} for a complete proof of this result.

The next result is an immediate corollary of Theorems \ref{thrm:discrete-conv} and \ref{thrm:nonconvergence-DT-DSGD}.
\begin{corollary}
Suppose 
the hypotheses of Theorem \ref{thrm:nonconvergence-DT-DSGD} are satisfied. Furthermore, suppose that every saddle point of $F$ is regular and satisfies Assumption \ref{a:eigvec-continuity}. Then, with probability 1, agents achieve asymptotic consensus and $x_n(k)$ converges to the set of local minima of $F$ for each agent $n$.
\end{corollary}

\section{Distributed Annealing Methods}
\label{sec:annealing-algorithms}

We will now consider first-order annealing algorithms for global optimization in $\R^d$. The distributed annealing algorithm we consider here is a simple generalization of D-SGD \eqref{eq:S-DGD}. At each iteration, each agent simply adds annealing noise to \eqref{eq:S-DGD}. If the annealing noise cools at an appropriate rate, the algorithm converges to a global minimum of \eqref{eq:f-distributed}.



We remark that the annealing algorithm discussed here is closely related to stochastic gradient Langevin dynamics (SGLD), popular for nonconvex machine learning tasks \cite{zhang2017hitting,raginsky2017non,chen2019}. The SGLD algorithm is a discretization of the Langevin diffusion \cite{chiang1987diffusion}. Recent work has studied the constant-weight variant of SGLD and characterized non-asymptotic properties of the algorithm such as local-minimum hitting and recurrence times \cite{zhang2017hitting,raginsky2017non,chen2019}.
Similar properties for the decaying weight parameter variant of SGLD were studied in \cite{chen2019}. When decaying weight parameters are used then SGLD coincides with the annealing algorithms considered in this paper.
Asymptotic properties of such algorithms including convergence to global minima and convergence rates were considered in \cite{gelfand1991recursive,yin1999rates}.

\subsection{Centralized SGD + Annealing} \label{sec:GD+annealing}
Before discussing the distributed annealing algorithm, it will be helpful to review centralized gradient-based annealing algorithms. In particular, we will see that if annealing noise is added to (centralized) SGD \eqref{eq:SGD-centralized}, then the process escapes from local minima and seeks out global minima. 

Given a differentiable function $f$, the centralized SGD + annealing algorithm is given by the recursion
\begin{equation} \label{eq:GD-plus-annealing}
x_{k+1} = x_k - \alpha_k( \nabla f(x_k) + \xi_k) + \gamma_{k}w_{k},
\end{equation}

\noindent where $\gamma_k$ denotes the annealing schedule and $w_k$ is i.i.d. $d$-dimensional Gaussian noise. 
Note that if the annealing noise is ``turned off'' (i.e., $\gamma_k$ is set to zero) then \eqref{eq:GD-plus-annealing} simply becomes SGD \eqref{eq:SGD-centralized}, and so converges to local minima. 

In order to ensure the convergence of SGD + annealing, additional structure must be assumed on the function $f$.
\begin{assumption}
\label{ass:GM_1} $f:\R^d\to\R$ is a $C^2$ function such that (i)  $\min_x f(x)$ exists, (ii)  $f(x)\to\infty$ and $\|\nabla f(x)\|\to\infty$ as $\|x\|\to\infty$, and (iii) $\inf_x (\|\nabla f(x)\|^2 - \Delta f(x) ) > -\infty$.
\end{assumption}
Part (i) of the assumption ensures that a global minimum of the problem  exists, while parts (ii)--(iii) are technical assumptions that ensure that, as the algorithm runs, the probability distribution corresponding to the position of $x_k$ cannot have mass escape out to infinity.

Our next assumption is a technical assumption regarding the structure of $f$. After presenting the assumption, we will discuss simple conditions under which it can be satisfied.
\begin{assumption}
\label{ass:GM_2} For $\e>0$ let $d\pi^\e(x) = \frac{1}{Z^\e}\exp\left(-\frac{2f(x)}{\e^2} \right)\dx,$
where
$
Z^\e= \int\exp\left(-\frac{2f(x)}{\e^2} \right)\dx,
$
and where $d\pi^e$ denotes the Radon-Nikodym derivative of a measure $\pi^\e$ taken with respect to the Lebesgue measure.
Assume $f$ is such that $\pi^\e$ has a weak limit $\pi$ as $\e\to 0$.
\end{assumption}
In the above assumption, the limiting distribution $\pi$ is constructed so as to place mass 1 on the set of global minima of $f$. The intuition underlying the assumption is that for fixed noise strength $\gamma_k =\e$, the distribution of $x_k$ can be shown to converge towards the stationary distribution $\pi^\e$ defined above \cite{chiang1987diffusion,gelfand1991recursive}. When sending $\gamma_k\to 0$, we would like the limiting distribution to be well defined and concentrate on the set of global minima---this is accomplished by Assumption \ref{ass:GM_2}.

The following lemma, adapted from \cite{hwang1980laplace}, Theorem 3.1, gives a simple condition on the Hessian of $f$ under which Assumption \ref{ass:GM_2} can be guaranteed to hold. In the lemma, $\lambda$ refers to the Lebesgue measure. 
\begin{lemma} \label{lemma:Laplace-ref}
Let
$\calN := \{x:f(x) = \inf_x f(x)\}.$
Suppose that (i) $\lambda(\{f(x)< a \}) > 0$ for any $a> \inf_x f(x)$, (ii) $\min_x f(x)$ exists and equals zero, (iii) There exists $\e>0$ such that $\{f(x) \leq \e\}$ is compact, (iv) $f$ is $C^3$.
Assume that $\calN$ consists of a finite set of isolated points and that the Hessian $\nabla^2 f(x)$ is invertible for all $x\in \calN$. Then the limit $\pi$ in Assumption~\ref{ass:GM_3} exists.
\end{lemma}

The next assumption imposes some asymptotic regularity on the objective $f$. 
\begin{assumption} \label{ass:GM_3} The following hold: 
\begin{itemize}
 \item [(i)] $\liminf_{|x|\to\infty}\langle \frac{\nabla f(x)}{\|\nabla f(x)\|}, \frac{x}{\|x\|} \rangle \geq C(d)$,
$C(d) = \left( \frac{4d-4}{4d-3} \right)^{\frac{1}{2}}$,\\
\item [(ii)] $\liminf_{\|x\|\to\infty} \frac{\|\nabla f(x)\|}{\|x\|} > 0$, 
\item [(iii)] $\limsup_{\|x\|\to\infty} \frac{\|\nabla f(x)\|}{\|x\|} < \infty$.
\end{itemize}
\end{assumption}
We emphasize Assumption \ref{ass:GM_3} holds asymptotically as $\|x\|\to \infty$. In applications, if a global minimum is known a priori to lie in some compact set, then this assumption can be trivially satisfied by modifying the objective function outside of the set.
See, for example \cite{swenson2019CAMSAP}.

We will assume that the noise process $\{\xi_k\}_{k\geq 1}$ satisfies Assumption \ref{a:noise-zero-mean} and that the annealing noise satisfies the following assumption, where we recall that $\calF_k$ represents the information available at iteration $k$.\footnote{Convergence under a relaxation of Assumption \ref{a:noise-zero-mean} is considered in \cite{gelfand1991recursive}. A similar relaxation in the distributed setting is considered in \cite{swenson2019CDC}}
\begin{assumption} \label{a:gauss}
$w_k$ is normally distributed, $w_k\sim \calN(0,I_d)$, independent of $\calF_k$.
\end{assumption}
Finally, we will assume that the weight sequences (and, in particular, the annealing schedule $\gamma_k$) take the following form.
\begin{assumption} \label{a:weights-annealing-centralized}
  $\alpha_k = c_\alpha \frac{1}{k}$ and $\gamma_k = c_\gamma (k\log{\log{k}})^{-\frac{1}{2}}$,
  where $c_\alpha$ and $c_\gamma$ are positive constants.
\end{assumption}

SGD + annealing converges in probability to the set of global minima. This is formalized in the following theorem from \cite{gelfand1991recursive}.
\begin{theorem}
\label{th:GM} Let $\{x_k\}_{k\geq 1}$ satisfy~\eqref{eq:GD-plus-annealing}. Suppose that Assumptions \ref{a:noise-zero-mean} and \ref{ass:GM_1}--\ref{a:weights-annealing-centralized} hold 
and assume $c_{\alpha}$ and $c_\gamma$ in Assumption \ref{a:weights-annealing-centralized} satisfy $c_\gamma/c_\alpha>C_{0}$, where the constant $C_0$ is defined after (2.3) in \cite{gelfand1991recursive}. Then $x_k$ converges in probability to the set of global minima of $f$ in the sense that, for any initial condition $x_0$, and for any $\e>0$ there holds $\mathbb{P}(\textup{dist}(x_k,S) > \e) \to 0$ as $k\to\infty$, where  $S = \arg\min f(x)$.
\end{theorem}

\subsection{Distributed Gradient Descent + Annealing} \label{sec:DGD+annealing}
We now consider distributed annealing algorithms for global optimization. 
Let $x_n(k)$ denote agent $n$'s estimate of an optimizer of \eqref{eq:f-distributed} at iteration $k$.
The D-SGD + annealing algorithm is defined agentwise by the recursion
\begin{equation} \label{eq:DGD_plus_annealing}
x_n(k+1) = x_n(k) + \beta_k\sum_{\ell\in \Omega_n} (x_\ell(k) - x_n(k)) - \alpha_k( \nabla f_n(x_n(k)) + \xi_n(k)) + \gamma_k w_n(k),
\end{equation}
for $n=1,\ldots,N$, where $w_n(k)$ is Gaussian noise injected independently at each agent and each time step, and $\gamma_k$ denotes the annealing schedule. If $\gamma_k$ is set to zero, then \eqref{eq:DGD_plus_annealing} becomes \eqref{eq:S-DGD} and the algorithm converges to local minima of \eqref{eq:f-distributed}.

In the distributed setting we will retain the assumptions from the centralized setting (but now applied to \eqref{eq:f-distributed}) and add the following assumption which ensures that agents can reach consensus.
\begin{assumption} \label{a:f_n-alignment}
For each $n$, the following hold: (i) $\nabla f_n$ is globally Lipschitz continuous, and (ii) there exists a $C>0$ such that $\langle x,\nabla f_n(x) \rangle \geq 0$ for all $\|x\| \geq C$.
\end{assumption}

We will also assume that the consensus weight sequence $\{\beta_k\}_{k\geq 1}$ takes the following form.
\begin{assumption}
\label{ass:weights-beta} $\beta_{k}=\frac{c_{\beta}}{k^{\tau_{\beta}}}$ for $k$ large, where $c_{\beta}>0$ and $\tau_{\beta}\in [0,1/2)$.
\end{assumption}

Under the above assumptions, agents reach consensus and converge to the set of \emph{global} minima of $F$, as stated next.
A proof of this result can be found in \cite{swenson2019CDC,swenson2020ICASSP}.
\begin{theorem}
\label{th:DGD-annealing1}
Let $\{(x_n(k))_{n=1}^N\}_{k\geq 1}$ satisfy \eqref{eq:DGD_plus_annealing} with arbitrary initial condition $(x_n(k))_{n=1}^N = x_0\in \R^{Nd}$. Suppose the communication graph satisfies Assumption \ref{a:G-connected-undirected}. Suppose also that Assumptions \ref{ass:GM_1}--\ref{ass:GM_3}  are satisfied with $f = \F$, where $\F$ is defined in \eqref{eq:f-distributed}, and that Assumptions \ref{a:f_n-alignment}--\ref{ass:weights-beta} hold. 
Suppose also that
for each $n$, the noise sequences $\{\xi_n(k)\}_{k\geq 1}$ and $\{w_n(k)\}_{k\geq 1}$ are mutually independent for all agents and satisfy Assumptions \ref{a:noise-zero-mean} and \ref{a:gauss} respectively, with $\xi_n(k) = \xi_k$ and $w_n(k) = w_k$.
Further, suppose that Assumption~\ref{a:weights-annealing-centralized} holds with 
$c_{\gamma}^{2}/c_{\alpha}>C_{0}$, where $C_0$ is defined after (2.3) in \cite{gelfand1991recursive}.
Then, $\lim_{k\to\infty}\|x_n(k) - x_\ell(k)\|=0$ with probability 1 for each $n,\ell$. Moreover, for each agent $n$, $\{x_n(k)_{k\geq 1}$ converges in probability to the set of global minima of $F$, i.e., for any $\e>0$, $\P(\textup{dist}(x_k,S)>\e) \to 0$ as $k\to\infty$, where $S = \arg\min F(x)$.
\end{theorem}

\section{Illustrative Examples} \label{sec:sims}
We illustrate the D-SGD and D-SGD + annealing algorithms on two simple examples: a small image classifier and a simple linear regression problem.

The image classifier is trained on the Fashion-MNIST task \cite{xiao2017fashion}.
This example simulates (on a smaller scale) the potential use-case of training large deep neural networks that require multiple machines to learn and where the training bottleneck is the communication bandwidth between machines.\footnote{The implementation is based on the distributed MNIST code from the examples section of the official PyTorch repository, https://github.com/pytorch/examples/tree/master/mnist\_hogwild} The model is a convolutional neural network with two convolutional layers followed by two fully connected layers and ReLU activation functions. All agents have their own independently shuffled copy of the data and the objective $f_n$ at agent $n$ is given by the empirical risk using the agent's data and the cross-entropy loss. 
Using the D-SGD algorithm, we train with batches of size 1000, a learning rate (parameter $\alpha_k$) of 0.5 decaying by 70\% each epoch, and the consensus term is set to $\beta_k=0.1$. Each node is allowed to make 10 full passes through the data. For the communication graph we consider both a cycle on four vertices and a 3-regular graph with 8 nodes. The results are shown in Table \ref{tab:fashion}, where we observe that the additional nodes improve the accuracy on the test set demonstrating the usefulness of the distributed algorithm. As a further comparison, we train a second centralized model for four times as many epochs so that it is a fair comparison with the decentralized run on four nodes. In this case the decentralized run is still able to beat the centralized baseline. It is important to remark that not all of the assumptions for D-SGD discussed in the paper are satisfied by this example (in particular, because we use ReLU activation functions, the function is not everywhere smooth).\footnote{Nonsmoothness, e.g., due to ReLU activation functions, is treated in \cite{swenson2020saddles}. Under mild assumptions D-SGD is known to converge to critical points. Saddle points are avoided under local regularity conditions near the saddle.} 
However, we find that that algorithm still performs well. An important future research direction is characterizing properties of distributed stochastic gradient algorithms under more general assumptions. 

\begin{table}[]
\centering
\begin{tabular}{cccc}
\textbf{Nodes} & \textbf{Epochs} & \textbf{Log. Loss} & \textbf{Accuracy} \\ \hline
1              & 10              & .484               & 81.7              \\
1              & 40              & .343               & 87.4              \\
4              & 10              & .318               & 88.0              \\
8              & 10              & .308               & 88.3
\end{tabular}
\caption{Performance of the D-SGD algorithm on the Fashion-MNIST dataset.}
\label{tab:fashion}
\end{table}

\begin{figure}[h]
    \centering
    \includegraphics[width=0.7\textwidth]{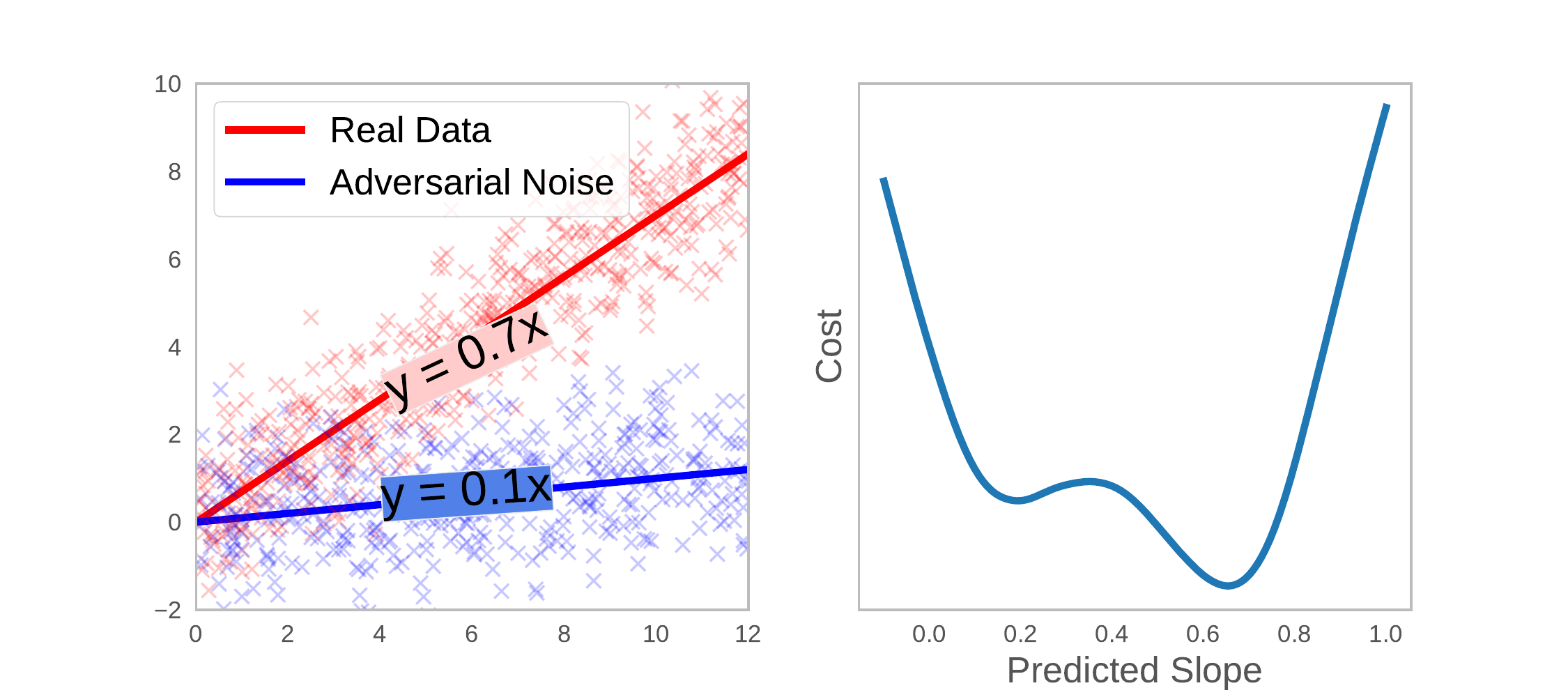}
    \caption{Linear regression data (left) and the associated optimization landscape (right).}
    \label{fig:lineplot}
\end{figure}

We now turn to the second example, an illustration of the D-SGD + annealing algorithm on a simple regression task. The linear regression model is uni-variate with no intercept ($y = w x$) with a non-convex loss function $L(y, \hat y) = \log(8(y - \hat y)^2 + 1)$ designed to give robustness to extreme outliers \cite{barron2019general}. Synthetic data is generated by sampling $x$ values uniformly in the interval $(0, 12)$ and using a binomial distribution with $p=0.55$ to pick between the ``real data" $y_1 = 0.7x + \epsilon$ and an adversarial distractor $y_2 = 0.1x + \epsilon$ where $\epsilon \sim \calN(0, 1)$. This distribution is contrived to have a sub-optimal local minimum near $w=0.1$ and a global minimum near $w=0.7$. The data and the loss curve are shown in Figure \ref{fig:lineplot}.
The objective at each agent $n$ is given by $f_n(w) = \frac{1}{N}\E_{{x,y}\sim\calD}(L(w x,y))$ and the overall objective is given by $F(w) = \E_{{x,y}\sim\calD}(L(wx,y))$, where $\calD$ represents the distribution for drawing samples previously mentioned. In this example, we consider an online learning setup. Let $(x_n(k), y_n(k))$ denote the sample drawn by agent $n$ in iteration $k$. The update step \eqref{eq:S-DGD} is done by taking $\nabla f_n(w_n(k)) + \xi_n(k) = \nabla \frac{1}{N}L(w_n(k) x_n(k),y_n(k))$.

We use two graphs: a cycle on four vertices and the Petersen graph, a 3-regular graph with 10 vertices and 15 edges. To implement D-SGD + annealing on this setup we use an exponentially decaying learning rate $\alpha_k = 0.01(0.998^k)$ and fix the consensus term as\footnote{While the theory holds for decaying $\beta_k$, it can be extended to constant $\beta_k=\beta$, for $\beta$ sufficiently small.} $\beta_k=4$ and annealing parameter $\gamma_k = 20 (0.9^{\sqrt{k}} )$. We found that, at least for this example, a more aggressive annealing schedule can be used than what is given in Assumption \ref{a:weights-annealing-centralized} in order to speed up convergence. We compare this to the case when there is no annealing term, $\gamma_k = 0$ (i.e., D-SGD). The simulation is run 100 times for $k=5000$ steps, and then we check to see which of the two minima it has converged to. In the case of the 4-Cycle graph, D-SGD with annealing converges to the global minimum in 90 out of 100 cases compared to 57 out of 100 when no annealing is used. For the Petersen graph, the global minimum is reached in 98 of 100 cases with annealing and 61 of 100 without.

\section{Conclusions} \label{sec:conclusion}
The paper reviewed refined convergence guarantees for distributed stochastic gradient algorithms. Convergence of D-SGD to local minima was discussed. The key role of the stable-manifold theorem in studying saddle-point nonconvergence was explored, along with the use of stochastic approximation techniques and continuous-time analysis methods. 
In order to obtain global optimality guarantees, D-SGD with additive annealing noise was discussed. 
The topic of refined convergence guarantees for distributed algorithms is a relatively new area and there are many open avenues for future research. We highlight a few now.  
Convergence rates for D-SGD + annealing are not currently understood, though convergence rate estimates do exist in the centralized case \cite{yin1999rates}.
If $\gamma_k$ in D-SGD + annealing is held constant, then the D-SGD + annealing algorithm is a distributed variant of SGLD. Nonasymptotic properties of these algorithms have been studied in the centralized setting \cite{zhang2017hitting,raginsky2017non,chen2019},
but not in the distributed setting.

\bibliographystyle{IEEEtran}
\bibliography{myRefs}

\end{document}